\documentclass[12pt]{article}
\usepackage{pstricks,pst-node}
\usepackage{psfig}
\usepackage{epsfig}
\usepackage{amsmath,amsthm}
\usepackage{amssymb}
\usepackage{color}
\usepackage{xspace}
\usepackage[colorlinks=true,
linkcolor=green,
filecolor=brown,
citecolor=green]{hyperref}

\def\gl{ground level\xspace}
\def\gf{generating function\xspace}
\catcode`\@=11

\thicklines
\newskip\Einheit \Einheit=.6cm
\newcount\xcoord \newcount\ycoord
\newdimen\xdim \newdimen\ydim \newdimen\PfadD@cke \newdimen\Pfadd@cke
\def\PfadDicke#1{\PfadD@cke#1 \divide\PfadD@cke by2 
\Pfadd@cke\PfadD@cke \multiply\PfadD@cke by2}
\long\def\LOOP#1\REPEAT{\def\BODY{#1}\ITERATE}
\def\ITERATE{\BODY \let\next\ITERATE \else\let\next\relax\fi \next}
\let\REPEAT=\fi
\def\Punkt{\hbox{\raise-2pt\hbox to0pt{\hss\scriptsize$\bullet$\hss}}}

\def\DuennPunkt(#1,#2){\unskip
  \raise#2 \Einheit\hbox to0pt{\hskip#1 \Einheit
          \raise-1.5pt\hbox to0pt{\hss\tiny$\bullet$\hss}\hss}}
		  
\def\NormalPunkt(#1,#2){\unskip
  \raise#2 \Einheit\hbox to0pt{\hskip#1 \Einheit
          \raise-3pt\hbox to0pt{\hss\large$\bullet$\hss}\hss}}
\def\DickPunkt(#1,#2){\unskip
  \raise#2 \Einheit\hbox to0pt{\hskip#1 \Einheit
          \raise-4pt\hbox to0pt{\hss\Large$\bullet$\hss}\hss}}
\def\Kreis(#1,#2){\unskip
  \raise#2 \Einheit\hbox to0pt{\hskip#1 \Einheit
          \raise-4pt\hbox to0pt{\hss\Large$\circ$\hss}\hss}}
\def\Diagonale(#1,#2)#3{\unskip\leavevmode
  \xcoord#1\relax \ycoord#2\relax
      \raise\ycoord \Einheit\hbox to0pt{\hskip\xcoord \Einheit
         \unitlength\Einheit
         \line(1,1){#3}\hss}}
\def\AntiDiagonale(#1,#2)#3{\unskip\leavevmode
  \xcoord#1\relax \ycoord#2\relax \advance\xcoord by -0.05\relax
      \raise\ycoord \Einheit\hbox to0pt{\hskip\xcoord \Einheit
         \unitlength\Einheit
         \line(1,-1){#3}\hss}}
\def\Pfad(#1,#2),#3\endPfad{\unskip\leavevmode
  \xcoord#1 \ycoord#2 \thicklines\ZeichnePfad#3\endPfad\thinlines}
\def\ZeichnePfad#1{\ifx#1\endPfad\let\next\relax
  \else\let\next\ZeichnePfad
    \ifnum#1=1
      \raise\ycoord \Einheit\hbox to0pt{\hskip\xcoord \Einheit
         \vrule height\Pfadd@cke width1 \Einheit depth\Pfadd@cke\hss}%
      \advance\xcoord by 1
     \else\ifnum#1=2
      \raise\ycoord \Einheit\hbox to0pt{\hskip\xcoord \Einheit
         \unitlength\Einheit
         \line(0,1){1}\hss}
      \advance\xcoord by 0
      \advance\ycoord by 1
 \else\ifnum#1=3
      \raise\ycoord \Einheit\hbox to0pt{\hskip\xcoord \Einheit
         \unitlength\Einheit
         \line(1,1){1}\hss}
      \advance\xcoord by 1
      \advance\ycoord by 1
    \else\ifnum#1=4
      \raise\ycoord \Einheit\hbox to0pt{\hskip\xcoord \Einheit
         \unitlength\Einheit
         \line(1,-1){1}\hss}
      \advance\xcoord by 1
      \advance\ycoord by -1
   \else\ifnum#1=5
      \raise\ycoord \Einheit\hbox to0pt{\hskip\xcoord \Einheit
         \unitlength\Einheit
         \line(2,1){2}\hss}
      \advance\xcoord by 2
      \advance\ycoord by 1
	  \else\ifnum#1=6
      \raise\ycoord \Einheit\hbox to0pt{\hskip\xcoord \Einheit
         \unitlength\Einheit
         \line(2,-1){2}\hss}
      \advance\xcoord by 2
      \advance\ycoord by -1
	  \else\ifnum#1=7
      \raise\ycoord \Einheit\hbox to0pt{\hskip\xcoord \Einheit
         \unitlength\Einheit
         \line(3,1){3}\hss}
      \advance\xcoord by 3
      \advance\ycoord by 1
	  \else\ifnum#1=8
      \raise\ycoord \Einheit\hbox to0pt{\hskip\xcoord \Einheit
         \unitlength\Einheit
         \line(3,-1){3}\hss}
      \advance\xcoord by 3
      \advance\ycoord by -1
    \fi\fi\fi\fi\fi\fi\fi\fi
  \fi\next}
\def\hSSchritt{\leavevmode\raise-.4pt\hbox 
to0pt{\hss.\hss}\hskip.2\Einheit
  \raise-.4pt\hbox to0pt{\hss.\hss}\hskip.2\Einheit
  \raise-.4pt\hbox to0pt{\hss.\hss}\hskip.2\Einheit
  \raise-.4pt\hbox to0pt{\hss.\hss}\hskip.2\Einheit
  \raise-.4pt\hbox to0pt{\hss.\hss}\hskip.2\Einheit}
\def\vSSchritt{\vbox{\baselineskip.2\Einheit\lineskiplimit0pt
\hbox{.}\hbox{.}\hbox{.}\hbox{.}\hbox{.}}}
\def\DSSchritt{\leavevmode\raise-.4pt\hbox to0pt{%
  \hbox to0pt{\hss.\hss}\hskip.2\Einheit
  \raise.2\Einheit\hbox to0pt{\hss.\hss}\hskip.2\Einheit
  \raise.4\Einheit\hbox to0pt{\hss.\hss}\hskip.2\Einheit
  \raise.6\Einheit\hbox to0pt{\hss.\hss}\hskip.2\Einheit
  \raise.8\Einheit\hbox to0pt{\hss.\hss}\hss}}
\def\dSSchritt{\leavevmode\raise-.4pt\hbox to0pt{%
  \hbox to0pt{\hss.\hss}\hskip.2\Einheit
  \raise-.2\Einheit\hbox to0pt{\hss.\hss}\hskip.2\Einheit
  \raise-.4\Einheit\hbox to0pt{\hss.\hss}\hskip.2\Einheit
  \raise-.6\Einheit\hbox to0pt{\hss.\hss}\hskip.2\Einheit
  \raise-.8\Einheit\hbox to0pt{\hss.\hss}\hss}}
\def\SPfad(#1,#2),#3\endSPfad{\unskip\leavevmode
  \xcoord#1 \ycoord#2 \ZeichneSPfad#3\endSPfad}
\def\ZeichneSPfad#1{\ifx#1\endSPfad\let\next\relax
  \else\let\next\ZeichneSPfad
    \ifnum#1=1
      \raise\ycoord \Einheit\hbox to0pt{\hskip\xcoord \Einheit
         \hSSchritt\hss}%
      \advance\xcoord by 1
    \else\ifnum#1=2
      \raise\ycoord \Einheit\hbox to0pt{\hskip\xcoord \Einheit
        \hbox{\hskip-2pt \vSSchritt}\hss}%
      \advance\ycoord by 1
    \else\ifnum#1=3
      \raise\ycoord \Einheit\hbox to0pt{\hskip\xcoord \Einheit
         \DSSchritt\hss}
      \advance\xcoord by 1
      \advance\ycoord by 1
    \else\ifnum#1=4
      \raise\ycoord \Einheit\hbox to0pt{\hskip\xcoord \Einheit
         \dSSchritt\hss}
      \advance\xcoord by 1
      \advance\ycoord by -1
    \fi\fi\fi\fi
  \fi\next}
\def\Koordinatenachsen(#1,#2){\unskip
 \hbox to0pt{\hskip-.5pt\vrule height#2 \Einheit width.5pt depth1 
\Einheit}%
 \hbox to0pt{\hskip-1 \Einheit \xcoord#1 \advance\xcoord by1
    \vrule height0.25pt width\xcoord \Einheit depth0.25pt\hss}}
\def\Koordinatenachsen(#1,#2)(#3,#4){\unskip
 \hbox to0pt{\hskip-.5pt \ycoord-#4 \advance\ycoord by1
    \vrule height#2 \Einheit width.5pt depth\ycoord \Einheit}%
 \hbox to0pt{\hskip-1 \Einheit \hskip#3\Einheit 
    \xcoord#1 \advance\xcoord by1 \advance\xcoord by-#3 
    \vrule height0.25pt width\xcoord \Einheit depth0.25pt\hss}}
\def\Gitter(#1,#2){\unskip \xcoord0 \ycoord0 \leavevmode
  \LOOP\ifnum\ycoord<#2
    \loop\ifnum\xcoord<#1
      \raise\ycoord \Einheit\hbox to0pt{\hskip\xcoord 
\Einheit\Punkt\hss}%
      \advance\xcoord by1
    \repeat
    \xcoord0
    \advance\ycoord by1
  \REPEAT}
\def\Gitter(#1,#2)(#3,#4){\unskip \xcoord#3 \ycoord#4 \leavevmode
  \LOOP\ifnum\ycoord<#2
    \loop\ifnum\xcoord<#1
      \raise\ycoord \Einheit\hbox to0pt{\hskip\xcoord 
\Einheit\Punkt\hss}%
      \advance\xcoord by1
    \repeat
    \xcoord#3
    \advance\ycoord by1
  \REPEAT}
\def\Label#1#2(#3,#4){\unskip \xdim#3 \Einheit \ydim#4 \Einheit
  \def\lo{\advance\xdim by-.5 \Einheit \advance\ydim by.5 \Einheit}%
  \def\llo{\advance\xdim by-.25cm \advance\ydim by.5 \Einheit}%
  \def\loo{\advance\xdim by-.5 \Einheit \advance\ydim by.25cm}%
  \def\o{\advance\ydim by.25cm}%
  \def\ro{\advance\xdim by.5 \Einheit \advance\ydim by.5 \Einheit}%
  \def\rro{\advance\xdim by.25cm \advance\ydim by.5 \Einheit}%
  \def\roo{\advance\xdim by.5 \Einheit \advance\ydim by.25cm}%
  \def\l{\advance\xdim by-.30cm}%
  \def\r{\advance\xdim by.30cm}%
  \def\lu{\advance\xdim by-.5 \Einheit \advance\ydim by-.6 \Einheit}%
  \def\llu{\advance\xdim by-.25cm \advance\ydim by-.6 \Einheit}%
  \def\luu{\advance\xdim by-.5 \Einheit \advance\ydim by-.30cm}%
  \def\u{\advance\ydim by-.30cm}%
  \def\ru{\advance\xdim by.5 \Einheit \advance\ydim by-.6 \Einheit}%
  \def\rru{\advance\xdim by.25cm \advance\ydim by-.6 \Einheit}%
  \def\ruu{\advance\xdim by.5 \Einheit \advance\ydim by-.30cm}%
  #1\raise\ydim\hbox to0pt{\hskip\xdim
     \vbox to0pt{\vss\hbox to0pt{\hss$#2$\hss}\vss}\hss}%
}
\catcode`\@=12

\begin{document}
\newtheorem{theorem}{Theorem}
\newtheorem{defn}[theorem]{Definition}
\newtheorem{lemma}[theorem]{Lemma}
\newtheorem{prop}[theorem]{Proposition}
\newtheorem{cor}[theorem]{Corollary}
\begin{center}
{\Large
An Identity for the Central Binomial Coefficient                         \\ 
}

\vspace{10mm}
David Callan  \\
\noindent {\small Dept. of Statistics, 
University of Wisconsin-Madison,  Madison, WI \ 53706}  \\
{\bf callan@stat.wisc.edu} 

June 14, 2012
\end{center}

\begin{abstract}
We find the joint distribution of three simple statistics on lattice paths of n upsteps and n downsteps leading to a triple sum identity for the central binomial coefficient {2n}-choose-{n}. We explain why one of the constituent double sums counts the irreducible pairs of compositions considered by Bender et al., and we evaluate some of the other sums.
\end{abstract}

\noindent{\Large \textbf{1 \  Introduction}  }
A \emph{Grand-Dyck path} is a lattice path consisting of an equal number of upsteps $U=(1,1)$ and downsteps $D=(1,-1)$. The horizontal line joining the endpoints is called \emph{\gl}\!\!. The number of upsteps is the \emph{semilength} of a Grand-Dyck path, also known as its \emph{size}. The number of Grand-Dyck paths of size $n$ is obviously the central binomial coefficient $\binom{2n}{n}$---choose locations for the upsteps among the $2n$ steps. A \emph{Dyck path} is a Grand-Dyck path that never goes below \gl, and it is \emph{primitive} if it is nonempty and its only vertices at \gl are its endpoints. 
The vertices at \gl of a nonempty Grand-Dyck path split it into \emph{components}, each of which is a primitive Dyck path or an inverted primitive Dyck path. 
A \emph{peak} in a Grand-Dyck path is an occurrence of $UD$ and a \emph{low} peak is one that starts at ground level.  A low peak is, in particular, a component of the path.

\Einheit=0.5cm
\[
\Pfad(-10,3),44434333343334434443\endPfad
\SPfad(-10,3),11111111111111111111\endSPfad
\DuennPunkt(-10,3)
\DuennPunkt(-9,2)
\DuennPunkt(-8,1)
\DuennPunkt(-7,0)
\DuennPunkt(-6,1)
\DuennPunkt(-5,0)
\DuennPunkt(-4,1)
\DuennPunkt(-3,2)
\DuennPunkt(-2,3)
\DuennPunkt(-1,4)
\DuennPunkt(0,3)
\DuennPunkt(1,4)
\DuennPunkt(2,5)
\DuennPunkt(3,6)
\DuennPunkt(4,5)
\DuennPunkt(5,4)
\DuennPunkt(6,5)
\DuennPunkt(7,4)
\DuennPunkt(8,3)
\DuennPunkt(9,2)
\DuennPunkt(10,3)
\Label\u{ \textrm{{\small A Grand-Dyck path of semilength 10, with 1 low peak, 2 components }}}(0,-0.8)
\Label\u{ \textrm{{\small above \gl, and 4 components altogether}}}(0,-1.7)
\]

\vspace*{4mm}

\noindent The \gf for Dyck paths counted by size is well known to be 
\[
C(x)=\frac{1-\sqrt{1-4x}}{2x}.
\]

In Section 2, we find the 4-variable \gf for Grand-Dyck paths counted by size, number of low peaks, number of components above \gl, and total number of components, and in Section 3, we find a closed formula for the joint distribution of these four statistics. In Section 4, we 
observe that the irreducible pairs of compositions considered by Bender et al. \cite{bender04} are equinumerous with low-peak-free Grand-Dyck paths, and we give a bijective explanation.

\vspace{5mm}

\noindent{\Large \textbf{2 \  A \gf}  }
Let $F(x,y,z,w)$ denote the \gf for Grand-Dyck paths where $x$ marks size, $y$ marks number of low peaks, $z$ marks number of components above \gl, and $w$ marks total number of components. The first return to \gl partitions nonempty Grand-Dyck paths into the 3 classes illustrated below, where $A$ denotes a Dyck path, $\overline{A}$ an inverted Dyck path, and $B$ a Grand-Dyck path.

\Einheit=0.5cm
\[
\Pfad(-9,0),34\endPfad
\Pfad(-3,0),3\endPfad
\Pfad(0,1),4\endPfad
\Pfad(5,0),4\endPfad
\Pfad(7,-1),3\endPfad
\Label\o{B}(-6.3,-0.3)
\Label\o{A\!\ne\!\epsilon}(-1,1)
\Label\o{\overline{A}}(6.5,-2)
\Label\o{B}(8.7,-0.3)
\Label\o{B}(1.7,-0.3)
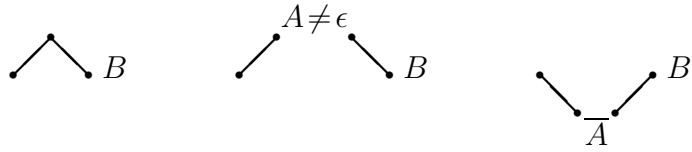
\Label\u{ \textrm{{\small A first return decomposition for nonempty Grand-Dyck paths}}}(0,-2.5)
\DuennPunkt(-9,0)
\DuennPunkt(-8,1)
\DuennPunkt(-7,0)
\DuennPunkt(-3,0)
\DuennPunkt(-2,1)
\DuennPunkt(0,1)
\DuennPunkt(1,0)
\DuennPunkt(6,-1)
\DuennPunkt(5,0)
\DuennPunkt(7,-1)
\DuennPunkt(8,0)
\]

\vspace*{5mm}

\noindent From this decomposition, we see that 
\[
F=1\ [\textrm{for the empty path}] +xyzwF +x(C(x)-1)zwF +xC(x)wF,
\]
an equation with solution
\begin{equation}\label{sol}
F(x,y,z,w)=\frac{2}{ 2 + 2 w x z (1-y) - w (1+z)(1- \sqrt{1-4 x})} \, .
\end{equation}
In particular, the \gf to count Grand-Dyck paths with no low peaks is 
\begin{equation}\label{nolowpeak}
F(x,0,1,1)=\frac{1}{x+\sqrt{1-4 x}} \, ,
\end{equation}
and the \gf to count Grand-Dyck paths by number of components above \gl is
\begin{equation}\label{compsAboveGL}
F(x,1,z,1)=\frac{2}{(z+1)\sqrt{1-4 x} -z+1}
\end{equation}
\vspace{5mm}

\noindent{\Large \textbf{3 \  An explicit count}  }
To obtain an explicit expression for the number $u(n,i,j,k)$ of Grand-Dyck paths of size $n$ with $i$ low peaks, $j$ components above \gl, and $k$ components altogether, first observe that there are $\binom{k}{i}$ ways to place the low peaks among the components. This reduces the problem to finding an expression for $v(n,j,k)$, the number of Grand-Dyck paths of size $n$ with no low peaks, $j$ components above \gl, and $k-j$ components below \gl. There are $\binom{k}{j}$ ways to arrange the above- and below-\gl components, so we may assume all components above \gl precede those below \gl. Each component above \gl has the form $UUPDQD$ where $P$ and $Q$ are (possibly empty) Dyck paths; each component below \gl has the form $D\overline{R}U$ where $\overline{R}$ is a Dyck path $R$ flipped over. Make the reversible transformations  $UUPDQD \rightarrow UPDUQD$ and $D\overline{R}\,U \rightarrow UR\,D$. Thus we see that the Grand-Dyck paths in question are equinumerous with Dyck paths of size $n$ and $2j+(k-j)=j+k$ components. It is well known that the number of Dyck paths of size $n$ with $k$ components is the generalized Catalan number $\frac{k}{2n-k}\binom{2n-k}{n-k}$ (arising as a $k$-fold convolution of Catalan numbers). Hence, 
\[
v(n,j,k)=\binom{k}{j} \frac{j+k}{2n-j-k}\binom{2n-j-k}{n-j-k},
\]
and, as noted above, $u(n,i,j,k)=\binom{k}{i}v(n-i,j-i,k-i)$, yielding
\[
u(n,i,j,k)=\binom{k}{i}\binom{k-i}{j-i}\frac{k-2i+j}{2n-j-k}\binom{2n-j-k}{n-i} \, .
\]
When $i=n$, as for the ``sawtooth'' path $(UD)^{n}$, we have $j=k=n$ and the expression for $u(n,i,j,k)$ is indeterminate; we must interpret it as 1.

Summing over $i,j,k$, we have the identity
\begin{equation}\label{id}
\binom{2n}{n}= \, 1\:+ \sum_{\substack{i,j,k \\[1mm] 0 \le i \le j \le k \le n \\[1mm] j+k < 2n}}\binom{k}{i}\binom{k-i}{j-i}\frac{k-2i+j}{2n-j-k}\binom{2n-j-k}{n-i} \, .
\end{equation}

\vspace{5mm}

{\Large \textbf{4 \  Irreducible pairs of compositions}  }
Bender et al. \cite{bender04} made the following definition. 

Let $n = b_{1} + \cdots + b_{k} = b_{1}' + \cdots + b_{k}'$ be a pair of compositions of $n$ into
$k$ positive parts. We say this pair is \emph{irreducible} if there is no positive $j < k$ for which
$b_{1} + \cdots + b_{j} = b_{1}'+ \cdots + b_{j}'$.

They showed that the number $f(n)$ of irreducible ordered pairs of compositions of $n$ into
the same (unspecified) number of parts has the \gf 
\begin{equation}\label{irred}
\sum_{n \geq 0} f(n+1) x^{n} ~=~ \frac{1}{x+\sqrt{1 - 4x}} \, .
\end{equation}
The generating functions in (\ref{nolowpeak}) and (\ref{irred}) are the same, which implies that irreducible pairs of compositions of $n+1$ are equinumerous with low-peak-free Grand-Dyck paths of size $n$.
It is not too hard to show this bijectively using a lattice path representation of compositions as illustrated below. Each entry $a_{i}$ in a composition $(a_{1}, \dots ,a_{k})$ contributes $a_{i}-1$ North steps followed by 1 East step. 
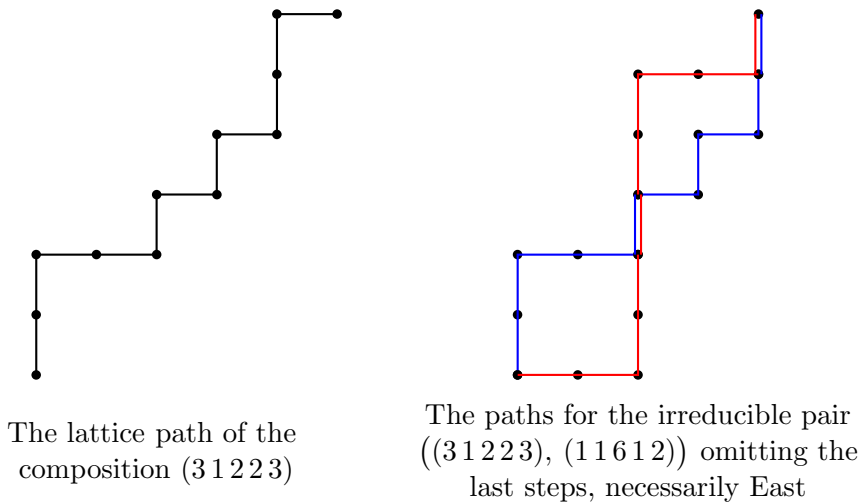
\begin{center}

\begin{pspicture}(-7.5,-1.5)(8.5,5)
\newrgbcolor{purple}{1 0 1}
 
\psset{unit=.8}   

\psdots(-6,0)(-6,1)(-6,2)(-5,2)(-4,2)(-4,3)(-3,3)(-3,4)(-2,4)(-2,5)(-2,6)(-1,6)
\psline(-6,0)(-6,1)(-6,2)(-5,2)(-4,2)(-4,3)(-3,3)(-3,4)(-2,4)(-2,5)(-2,6)(-1,6)

\psdots(2,0)(2,1)(2,2)(3,2)(4,2)(4,3)(5,3)(5,4)(6,4)(6,5)(6,6)
\psline[linecolor=blue](2,0)(2,1)(2,2)(3,2)(4,2)
\psline[linecolor=blue](4,3)(5,3)(5,4)(6,4)(6,5)

\psdots(2,0)(3,0)(4,0)(4,1)(4,2)(4,3)(4,4)(4,5)(5,5)(6,5)
\psline[linecolor=red](2,0)(3,0)(4,0)(4,1)(4,2)
\psline[linecolor=red](4,3)(4,4)(4,5)(5,5)(6,5)

\psline[linecolor=red](4.05,2)(4.05,3)
\psline[linecolor=blue](3.95,2)(3.95,3)

\psline[linecolor=blue](6.05,5)(6.05,6)
\psline[linecolor=red](5.95,5)(5.95,6)


\rput(-4,-1){\small The lattice path of the }
\rput(-4,-1.6){\small composition $(3\,1\,2\,2\,3)$}

\rput(4,-.7){\small The paths for the irreducible pair}
\rput(4,-1.3){\small $\big( (3\,1\,2\,2\,3),\,(1\,1\,6\,1\,2)\big)$ omitting the}
\rput(4,-1.9){\small last  steps, necessarily East}

\end{pspicture}

\end{center} 
The diagram on the right above represents an irreducible ordered pair of compositions. By definition of irreducible, no East step in the first (blue) path coincides with an East step in the second (red) path.
The vertices common to both paths split the diagram into path pairs that form parallelogram polyominoes, possibly the degenerate polyomino consisting of 2 coincident North steps. There are several bijections \cite[Ex. 6.19 $\ell$]{ec2} from parallelogram polyominoes of size (semiperimeter) $k$ to Dyck paths of size $k-1$. By elevating the resulting Dyck path (prepend an upstep, append a downstep), we get a size-preserving bijection from parallelogram polyominoes to primitive Dyck paths. So apply this bijection to each parallelogram polyomino to get a primitive Dyck path and use the color of the upper path to determine whether to flip it over. 
The degenerate polyomino corresponds to the Dyck path $UD$ and we always flip this over  because we don't want any low peaks. Lastly, concatenate the Dyck paths to obtain a low-peak-free Grand-Dyck path. 
This is the desired bijection from irreducible pairs of compositions of $n+1$ to low-peak-free Grand-Dyck paths of size $n$.

Emanuele Munarini notes in the OEIS \cite{oeis} entry  for the sequence 
\htmladdnormallink{A081696}{http://oeis.org/A081696} that the number of irreducible pairs of compositions of $n+1$
can be expressed as $\sum_{j=0}^{n}\frac{3j+1}{n+j+1}\binom{2n-j}{n-2j}$. In fact, our results can be used to get the following counts.
\begin{enumerate}
\item The number of low-peak-free Grand-Dyck paths with $j$ components above \gl is Munarini's summand $\frac{3j+1}{n+j+1}\binom{2n-j}{n-2j}$. 

\item The number of low-peak-free Grand-Dyck paths with $j$ components above \gl and $k$ components 
altogether is $\frac{j+k}{2n-j-k} \binom{2n-j-k}{n-j-k}\binom{k}{j}$. 

\item The number of unrestricted Grand-Dyck paths with $j$ components above \gl is $\frac{2j+1}{2n+1}\binom{2n+1}{n-j}$.

\item The number of unrestricted Grand-Dyck paths with $j$ components above \gl and $k$ components 
altogether is $\frac{k}{2n-k}\binom{2n-k}{n-k}\binom{k}{j}$.

\item The number of unrestricted Grand-Dyck paths with $j$ big components above \gl is $\frac{2j+1}{n+1}\binom{2n+2}{n-2j}$ (a big component is one of size $\ge 2$).
\end{enumerate}
Are there bijective proofs?

\end{document}